\documentclass[11pt]{amsart}
\usepackage{amsmath}
\usepackage{amssymb}
\usepackage{amscd}
\usepackage{color}
\usepackage[top=1in, bottom=1in, left=2.5cm, right=2.5cm]{geometry}

\def\NZQ{\mathbb}               
\def\NN{{\NZQ N}}
\def\QQ{{\NZQ Q}}
\def\ZZ{{\NZQ Z}}
\def\RR{{\NZQ R}}
\def\CC{{\NZQ C}}

\newtheorem{Theorem}{Theorem}[section]
\newtheorem{Lemma}[Theorem]{Lemma}
\newtheorem{Corollary}[Theorem]{Corollary}
\newtheorem{Proposition}[Theorem]{Proposition}

\newtheorem{Example}[Theorem]{Example}

\newtheorem{Definition}[Theorem]{Definition}

\let\epsilon\varepsilon
\let\phi=\varphi
\let\kappa=\varkappa

\begin{document}
\title{Positivity of mixed multiplicities of  filtrations}
\author{Steven Dale Cutkosky}

\author{Hema Srinivasan}
\author{Jugal Verma}

\thanks{This work was done while the first two authors were Visiting Professors at the Indian Institute of Technology, Bombay}
\thanks{The first author was partially supported by  NSF grant DMS-1700046}
\address{Steven Dale Cutkosky, Department of Mathematics,
University of Missouri, Columbia, MO 65211, USA}
\email{cutkoskys@missouri.edu}

\address{Hema Srinivasan, Department of Mathematics,
University of Missouri, Columbia, MO 65211, USA}
\email{srinivasanh@missouri.edu}

\address{Jugal Verma, Department of Mathematics, 
Indian Institute of Technology, Bombay, Mumbai 40076, India}
\email{jkv@math.iitb.ac.in}

\begin{abstract} The theory of mixed multiplicities of filtrations by $m$-primary ideals in a ring is introduced in \cite{CSS}.   In this paper, we consider the positivity of mixed multiplicities of filtrations.  We show that the mixed multiplicities of filtrations must be nonnegative real numbers  and give examples to show that they could be zero or even irrational.  When $R$ is analytically irreducible,  and   $\mathcal I(1),\ldots,\mathcal I(r)$ are filtrations of  $R$ by $m_R$-primary ideals,  we show that all of the mixed multiplicities $e_R(\mathcal I(1)^{[d_1]},\ldots,\mathcal I(r)^{[d_r]};R)$    are positive if and only if the ordinary  multiplicities  $e_R(\mathcal I(i);R)$ for $1\le i\le r$ are positive.  We extend this to modules and prove a simple characterization of when the mixed multiplicities are positive or zero on a finitely generated module. 
\end{abstract}


\maketitle
\section{Introduction}

The study of mixed multiplicities of $m_R$-primary ideals in a Noetherian local ring $R$ with maximal ideal $m_R$  was initiated by Bhattacharya \cite{Bh}, Rees  \cite{R} and Teissier  and Risler \cite{T1}. 
In \cite{CSS}  the notion of mixed multiplicities  is extended to arbitrary,  not necessarily Noetherian, filtrations of $R$ by $m_R$-primary ideals. It is shown in \cite{CSS} that many basic theorems for mixed multiplicities of $m_R$-primary ideals hold true for filtrations. 

The development of the subject of mixed multiplicities and its connection to Teissier's work on equisingularity \cite{T1} can be found in \cite{GGV}.   A  survey of the theory of  mixed multiplicities of  ideals  can be found in  \cite[Chapter 17]{HS}, including discussion of the results of  the papers \cite{R1} of Rees and \cite{S} of  Swanson, and the theory of Minkowski inequalities of Teissier \cite{T1}, \cite{T2}, Rees and Sharp \cite{RS} and Katz \cite{Ka}.   Later, Katz and Verma \cite{KV}, generalized mixed multiplicities to ideals which are not all $m_R$-primary.  Trung and Verma \cite{TV} computed mixed multiplicities of monomial ideals from mixed volumes of suitable polytopes.  Mixed multiplicities are also used by Huh in the analysis of the coefficients of the chromatic polynomial of graph theory in \cite{H}.

We will be concerned with  multiplicities and mixed multiplicities of  (not necessarily Noetherian) filtrations, which are defined as follows. 

\begin{Definition}
A filtration $\mathcal I=\{I_n\}_{n\in\NN}$ of  a ring $R$ is a descending chain
$$
R=I_0\supset I_1\supset I_2\supset \cdots
$$
of ideals such that $I_iI_j\subset I_{i+j}$ for all $i,j\in \NN$.  A filtration $\mathcal I=\{I_n\}$ of  a local ring $R$ by $m_R$-primary ideals is a filtration $\mathcal I=\{I_n\}_{n\in\NN}$ of $R$ such that   $I_n$ is $m_R$-primary for $n\ge 1$.
A filtration $\mathcal I=\{I_n\}_{n\in\NN}$ of  a ring $R$ is said to be Noetherian if $\bigoplus_{n\ge 0}I_n$ is a finitely generated $R$-algebra.
 \end{Definition}

The key result needed to define the multiplicity of a filtration of $R$ by $m_R$-primary ideals is the following.  Let $\ell_R(M)$ denote the length of an $R$-module $M$.

\begin{Theorem} \label{TheoremI20} (\cite[Theorem 1.1]{C2} and  \cite[Theorem 4.2]{C3}) Suppose that $R$ is a Noetherian local ring of dimension $d$, and  $N(\hat R)$ is the nilradical of the $m_R$-adic completion $\hat R$ of $R$.  Then   the limit 
\begin{equation}\label{I5}
\lim_{n\rightarrow\infty}\frac{\ell_R(R/I_n)}{n^d}
\end{equation}
exists for any filtration  $\mathcal I=\{I_n\}$ of $R$ by $m_R$-primary ideals, if and only if $\dim N(\hat R)<d$.
\end{Theorem}

When the ring $R$ is a domain and is essentially of finite type over an algebraically closed field $k$ with $R/m_R=k$, Lazarsfeld and Musta\c{t}\u{a} \cite{LM} showed that
the limit exists for all filtrations  of $R$ by $m_R$-primary ideals.  Cutkosky \cite{C3} proved it in the complete generality as stated above in Theorem \ref{TheoremI20}.

As can be seen from this theorem,  one must impose the condition that 
the dimension of the nilradical of the completion $\hat R$ of $R$ is less than the dimension of $R$. The nilradical $N(R)$ of a $d$-dimensional ring $R$ is 
$$
N(R)=\{x\in R\mid x^n=0 \mbox{ for some positive integer $n$}\}.
$$
We have that $\dim N(R)=d$ if and only if there exists a minimal prime $P$ of $R$ such that $\dim R/P =d$ and $R_P$ is not reduced. In particular,  the condition $\dim N(\hat R)<d$ holds if $R$ is analytically unramified; that is, $\hat R$ is reduced.

The multiplicity of a non Noetherian filtration can be an irrational number.  We will now give a
 very simple example of a filtration by $m_R$-primary ideals with an irrational multiplicity. 
 Let $k$ be a field and $R=k[[x]]$ be a power series ring over $k$. Let $I_n =(x^{\lceil n\sqrt{2}\rceil})$   
where $\lceil \alpha\rceil$ is the round up of a real number $\alpha$ (the smallest integer which is greater than or equal to $\alpha$).  Then $\{I_n\}$ is a graded family of $m_R$-primary ideals such that 
$$
\lim_{n\rightarrow\infty}\frac{\ell_R(R/I_n)}{n} =\sqrt{2}
$$ 
is an irrational number. 

Mixed multiplicities of filtrations are defined in \cite{CSS}. 
 Let $M$ be a finitely generated $R$-module where $R$ is a $d$-dimensional Noetherian local ring with $\dim N(\hat R)<d$. Let $\mathcal I(1)=\{I(1)_n\},\ldots, \mathcal I(r)=\{I(r)_n\}$ be filtrations of $R$ by $m_R$-primary ideals. 
 In  \cite[Theorem 6.1]{CSS} and  \cite[Theorem 6.6]{CSS}, it is shown that the function
\begin{equation}\label{M2}
P(n_1,\ldots,n_r)=\lim_{m\rightarrow \infty}\frac{\ell_R(M/I(1)_{mn_1}\cdots I(r)_{mn_r}M)}{m^d}
\end{equation}
is equal to a homogeneous polynomial $G(n_1,\ldots,n_r)$ of total degree $d$ with real coefficients for all  $n_1,\ldots,n_r\in\NN$.  

We  define the mixed multiplicities of $M$ from the coefficients of $G$, generalizing the definition of mixed multiplicities for $m_R$-primary ideals. Specifically,   
 we write 
$$
G(n_1,\ldots,n_r)=\sum_{d_1+\cdots +d_r=d}\frac{1}{d_1!\cdots d_r!}e_R(\mathcal I(1)^{[d_1]},\ldots, \mathcal I(r)^{[d_r]};M)n_1^{d_1}\cdots n_r^{d_r}.
$$
We say that $e_R(\mathcal I(1)^{[d_1]},\ldots,\mathcal I(r)^{[d_r]};M)$ is the mixed multiplicity of $M$ of type $(d_1,\ldots,d_r)$ with respect to the filtrations $\mathcal I(1),\ldots,\mathcal I(r)$.
Here we are using the notation 
\begin{equation}\label{eqI6}
e_R(\mathcal I(1)^{[d_1]},\ldots, \mathcal I(r)^{[d_r]};M)
\end{equation}
  to be consistent with the classical notation for mixed multiplicities of $M$ for $m_R$-primary ideals from \cite{T1}. The mixed multiplicity of $M$ of type $(d_1,\ldots,d_r)$ with respect to $m_R$-primary ideals $I_1,\ldots,I_r$, denoted by $e_R(I_1^{[d_1]},\ldots,I_r^{[d_r]};M)$ (\cite{T1}, \cite[Definition 17.4.3]{HS}) is equal to the mixed multiplicity $e_R(\mathcal I(1)^{[d_1]},\ldots,\mathcal I(r)^{[d_r]};M)$, where the Noetherian $I$-adic filtrations $\mathcal I(1),\ldots,\mathcal I(r)$ are defined by $\mathcal I(1)=\{I_1^i\}_{i\in \NN}, \ldots,\mathcal I(r)=\{I_r^i\}_{i\in \NN}$.

We write the multiplicity $e_R(\mathcal I;M)=e_R(\mathcal I^{[d]};M)$ if $r=1$, and $\mathcal I=\{I_i\}$ is a filtration of  $R$ by $m_R$-primary ideals. We have that
$$
e_R(\mathcal I;M)=\lim_{m\rightarrow \infty}d!\frac{\ell_R(M/I_mM)}{m^d}.
$$

 Valuation ideals give   natural examples of  filtrations.
  Suppose that $R$ is a $d$-dimensional  excellent  local domain. A valuation $\nu$ of the quotient field of $R$ is called divisorial if  the valuation ring  $V_{\nu}$ of $\nu$ dominates a localization of $R$ at a nonzero prime ideal $P$ of $R$ ($R_P\subset V_{\nu}$ and $m_{\nu}\cap R_P=P_P$) and $V_{\nu}$ is essentially of finite type over $R$ ($V_{\nu}$ is a localization of a finitely generated $R$-algebra). We have that $\nu$ is divisorial if and only if there exists a normal projective $R$-scheme $X$ with a birational projective morphism $\pi:X\rightarrow \mbox{Spec}(R)$ and a codimension one closed subvariety $E$ of $X$ such that the local ring $\mathcal O_{X,E}=V_{\nu}$ is the valuation ring of $\nu$.
 Define valuation ideals
$I(\nu)_n=\{f\in R\mid \nu(f)\ge n\}$ in $R$ for $n\in\NN$.

Suppose that $\nu$ is a divisorial valuation which dominates  $R$. Then $\nu$ determines a filtration $\mathcal I(\nu)$ of $R$ by $m_R$-primary ideals, by   
$\mathcal I(\nu)=\{I(\nu)_n\}$.
In  a two dimensional normal local ring $R$, the condition that the  filtration of valuation ideals in $R$ is Noetherian for all  divisorial valuations dominating $R$ is the condition (N) of Muhly and Sakuma \cite{MS}. It is proven in \cite{C4} that a complete normal local ring of dimension two satisfies condition (N) if and only if it's divisor class group is a torsion group. 
It follows from   \cite[Theorem 9]{CS} that the  multiplicity $e_R(\mathcal I(\nu);R)$ of the filtration of a 
divisorial valuation $\nu$ dominating a two dimensional excellent and normal local ring $R$    is always a rational  number.
However, in dimension three it can happen that the multiplicity of the filtration of a valuation can be irrational. In \cite[Example 6]{CS}, an example is given of a divisorial valuation $\nu$ dominating an excellent local domain $R$ of dimension three such that $e_R(\mathcal I(\nu);R)$ is an irrational number.

Suppose that $\nu_1,\ldots,\nu_r$ are divisorial valuations of the quotient field  of $R$ which dominate $R$. 
Then  for $n_1,\ldots,n_r\in \NN$, the function
\begin{equation}\label{eqR12}
G(n_1,\ldots,n_r)=\lim_{n\rightarrow\infty}\frac{\ell_R(R/I(\nu_1)_{nn_1}\cdots I(\nu_r)_{n n_r})}{n^d}
\end{equation}
of equation (\ref{M2}) is a homogeneous  polynomial of total degree $d$, whose coefficients determine the mixed multiplicities $e_R(\mathcal I(\nu_1)^{[d_1]},\ldots,\mathcal I(\nu_r)^{[d_r]};R)$ of (\ref{eqI6}).

It can be deduced from the rationality of the multiplicities $e_R(\mathcal I(\nu_i);R)$ in dimension two  that the mixed multiplicities of valuation ideals in a two dimensional excellent and normal local ring   are always rational  numbers; that is, the coefficients of (\ref{eqR12}) are always rational numbers if $R$ has dimension two.  However, the mixed multiplicities of valuation ideals  can be irrational if $d\ge 3$, since the multiplicities $e_R(\mathcal I(\nu_i);R)$ can be irrational.

 Using methods of Rees as in the proof of formula (8) of \cite{C1}, we can deduce that the mixed multiplicities $e_R(\mathcal I(\nu_1)^{[d_1]},\ldots,\mathcal I(\nu_r)^{[d_r]};R)$ are always positive if $\nu_1,\ldots,\nu_r$ are divisorial valuations which dominate an excellent analytically irreducible local domain. 


In the classical case of $m_R$-primary ideals, we also have that all mixed multiplicities are positive. 
If $R$ is a  $d$-dimensional Noetherian local ring, $I$ is an $m_R$-primary ideal in $R$ and $M$ is a finitely generated $R$-module of dimension $d$ then the multiplicity
$e_R(I;M)>0$. Further, if $J_1,\ldots,J_r$ are $m_R$-primary ideals,   then all mixed multiplicities $e_R( J_1^{[d_1]},\ldots, J_r^{[d_r]};M)$ are positive if $\dim M=d$ (\cite{T1} or \cite{HS}[Corollary 17.4.7]).

In contrast,  if $R$ is a $d$-dimensional Noetherian local ring such that $\dim N(\hat R)<d$ and $\mathcal I=\{I_n\}$ is a filtration of $m_R$-primary ideals,  then the limit 
$$
e_R(\mathcal I;R)=d!\lim_{n\rightarrow\infty}\frac{\ell_R(R/I_n)}{n^d}
$$
can be zero if the filtration is non Noetherian. A simple example is the filtration $\mathcal I=  \{I_n\}$ where $I_n=(x_1)+m_R^n$ in $R=\CC[[x_1,\ldots,x_d]]$ (with $d\ge 2$).

The mixed multiplicities of filtrations are always nonnegative, as we show in the following proposition. 

\begin{Proposition}\label{Prop0} Suppose that $R$ is a Noetherian local ring of dimension $d$ such that $\dim N(\hat R)<d$. Suppose that $\mathcal I(1),\ldots,\mathcal I(r)$ are filtrations of $R$ by $m_R$-primary ideals, and $M$ is a finitely generated $R$-module. Then 
 for all $d_1,\ldots,d_r$ with $d_1+\cdots+d_r=d,$
the mixed multiplicities $e_R(\mathcal I(1)^{[d_1]},\ldots,\mathcal I(r)^{[d_r]};M)$ are nonnegative real numbers.
\end{Proposition}

A natural question, at this point, is  whether the mixed multiplicities are always strictly positive if the multiplicities $e_R( \mathcal I(j);R)$ are positive . This is in fact true if $R$ is analytically irreducible, as we  show in the following theorem.

\begin{Theorem}\label{Theorem2} Suppose that $R$ is a $d$-dimensional  analytically irreducible Noetherian local ring and
$\mathcal I(1)=\{I(1)_n\},\ldots,\mathcal I(r)=\{I(r)_n\}$ are filtrations of $R$ by $m_R$-primary ideals such that 
$$
e_R(\mathcal I(j);R)=d!\lim_{n\rightarrow\infty}\frac{\ell_R(R/I(j)_n)}{n^d}>0
$$
for  $1\le j\le r$. Then all of the mixed multiplicities
$$
e_R(\mathcal I(1)^{[d_1]},\ldots,\mathcal I(r)^{[d_r]};R)
$$
 (for all $d_1,\ldots,d_r\in \NN$ such that $d_1+\cdots+d_r=d$)
 are  positive.
 \end{Theorem}
 
 However, there do exist excellent domains for which  all $e_R( \mathcal I(j);R)$ are positive but not all of the mixed multiplicities are positive. We give an example,  Example \ref{Example1}, which is established in Section \ref{Sec3}.

We have the following corollary to Theorem \ref{Theorem2}, giving general conditions for all mixed multiplicities of filtrations of $m_R$-primary ideals to be positive.
 
\begin{Corollary}\label{Cor3} Suppose that $R$ is a Noetherian local ring of dimension $d$ with $\dim N(\hat R)<d$. Suppose $\mathcal I(j)=\{I(j)_i\}$ for $1\le j\le r$ are filtrations of $R$ by $m_R$-primary ideals and  $M$ is a finite  $R$-module  of dimension $d$ and $\mathcal I(1),\ldots \mathcal I(r)\}$ are filtrations of $R$ by $m_R$-primary ideals.  Suppose that 
$$
e_{\hat R/P}(\mathcal I(j)\hat R/P;\hat R/P)>0
$$
for $1\le j\le r$ and all minimal primes $P$ of $\hat R$ such that $\dim \hat R/P=d$. Then all of the mixed multiplicities $$
e_R(\mathcal I(1)^{[d_1]},\ldots,\mathcal I(r)^{[d_r]};M)
$$
for all  $d_1,\ldots,d_r\in \NN$ such that $d_1+\cdots+d_r=d$
 are  positive.
\end{Corollary}

Proofs of the above results are given in Section \ref{Sec3}.

We generalize this to all analytically irreducible local rings.  We obtain the following necessary and sufficient criterion for vanishing and positivity of mixed multiplicities of filtrations. 

Given filtrations $\mathcal I(1),\ldots,\mathcal I(r)$ of $R$ by $m_R$-primary ideals, we can reindex them so that
there is an $s$ with $0\le s\le r$ such that 
$e_R(\mathcal I(j);R)>0$ for $1\le j\le s$ and $e_R(\mathcal I(j);R)=0$ for $s<j\le r$.

\begin{Theorem}\label{Theorem3} Suppose that $R$ is a $d$-dimensional  analytically irreducible Noetherian local ring, $M$ is a finitely generated $R$-module of dimension $d$ and
$\mathcal I(1)=\{I(1)_n\},\ldots,\mathcal I(r)=\{I(r)_n\}$ are filtrations of $m_R$-primary ideals such that there is 
an $s$ with $0\le s\le r$ such that 
$e_R(\mathcal I(j);R)>0$ for $1\le j\le s$ and $e_R(\mathcal I(j);R)=0$ for $s<j\le r$.
Then the mixed multiplicities
\begin{equation}\label{eqX41}
e_R(\mathcal I(1)^{[d_1]},\ldots,\mathcal I(r)^{[d_r]};M)=
0\mbox{ if }d_{s+1}+\cdot+d_r>0
\end{equation}
and
\begin{equation}\label{eqX5}
e_R(\mathcal I(1)^{[d_1]},\ldots,\mathcal I(r)^{[d_r]};M)=e_R(\mathcal I(1)^{[d_1]},\ldots,\mathcal I(s)^{[d_s]};M)> 0\mbox{ if }d_{s+1}+\cdots+d_r=0
\end{equation}
for all  $d_1,\ldots,d_r\in \NN$ such that $d_1+\cdots+d_r=d$.
 \end{Theorem}
 
 We have the following immediate corollary.
 
 \begin{Corollary}\label{Cor4} Suppose that $R$ is an analytically irreducible Noetherian local ring of dimension $d$, $M$ is a finite $R$-module of dimension $d$ and
$\mathcal I(1)=\{I(1)_n\},\ldots,\mathcal I(r)=\{I(r)_n\}$ are filtrations of $R$ by $m_R$-primary ideals such that $e_R(\mathcal I(j);R)=0$ for $1\le j\le r$.
Then the mixed multiplicities
$$
e_R(\mathcal I(1)^{[d_1]},\ldots,\mathcal I(r)^{[d_r]};M)=0
$$
for all  $d_1,\ldots,d_r\in \NN$ such that $d_1+\cdots+d_r=d$.
 \end{Corollary}
 In the case that $r=2$, Corollary \ref{Cor4} follows directly from the third Minkowski inequality for filtrations of \cite[Theorem 6.3]{CSS}.

 Theorem \ref{Theorem3} is proved in Section \ref{Sec5} of this paper. 
  
Throughout this paper,   $\NN$ will denote the non-negative integers and   $\ZZ_+$ will denote the positive integers.  We will denote the set of nonnegative rational numbers  by $\QQ_{\ge 0}$, the positive rational numbers by $\QQ_+$, and the set of non-negative real numbers by $\RR_{\ge0}$.

For a local ring $R$,   $m_R$ denotes the maximal ideal. The quotient field of a domain $R$ will be denoted by ${\rm QF}(R)$.

\section{Mixed multiplicities on complete local domains}\label{Sec2}
 Suppose that $R$ is a complete Noetherian local domain of dimension $d$, and $\mathcal I=\{I_n\}$ is a filtration of $R$ by  $m_R$-primary ideals.  
 
 For $a\in \ZZ_+$, 
 let $\mathcal I_a=\{I_{a,i}\}$ be the $a$-th truncated filtration of $\mathcal I$ defined in \cite[Definition 4.1]{CSS}. 

\begin{Definition}\label{trunc}
 Suppose that $\mathcal I=\{I_i\}$  is a filtration of a  local ring $R$. For $a\in \ZZ_+$, the $a$-th truncated  filtration $\mathcal I_a=\{I_{a,i}\}$ of $\mathcal I$ is defined  by $I_{a,n}=I_n$ if $n\le a$ and if $n>a$, then $I_{a,n}=\sum I_{a,i}I_{a,j}$ where the sum is over $i,j>0$ such that $i+j=n$.  
 \end{Definition}

 For $s\in \ZZ_+$, let $\mathcal I[s]$ denote the filtration $\mathcal I[s]=\{I_{si}\}$.

We first review a method for computing asymptotic multiplicities, developed in  \cite{C1}, \cite{C2}, \cite{C3} and \cite{CSS}.  The method is inspired by the work of \cite{Ok}, \cite{LM} and \cite{KK} on volumes of  linear series. There exists a regular local ring $S$ of dimension $d$ which is a localization of a finitely generated $R$-algebra with the same quotient field ${\rm QF}(R)$ as $R$, which dominates $R$ ($R\subset S$ and $m_S\cap R=m_R$). An algebraic  proof of this is given in \cite[Lemma 4.2]{CSS}.  Letting $y_1,\ldots,y_d$ be a regular system of parameters in $S$, we define a valuation $\nu$ dominating $S$ by prescribing that $\nu(y_i)=\lambda_i$ for $1\le i\le d$, where $\lambda_i\in \RR$ are linearly  independent over the  field $\QQ$ of rational numbers and satisfy $\lambda_i\ge 1$ for all $i$. Let $V_{\nu}$ be the valuation ring of $\nu$ and
for $\lambda\in \RR_{\ge 0}$, let
$$
K_{\lambda}=\{f\in {\rm QF}(R)\mid \nu(f)\ge \lambda\}
$$
and 
$$
K_{\lambda}^+=\{f\in {\rm QF}(R)\mid \nu(f)> \lambda\}
$$
which are ideals in $V_{\nu}$. Let $k=R/m_R$. 

There exists $\overline c\in \ZZ_+$ such that $m_R^{\overline c}\subset I_1$, so that $m_R^{n\overline c}\subset I_n$ for all $n$. By equation (10) of \cite{C1} or equation (31) of \cite{C3}, there exists $\beta\in \ZZ_+$ such that 
\begin{equation}\label{eq1}
K_{\beta n}\cap R\subset m_R^{n\overline c}\subset I_n
\end{equation}
for all $n$.

\begin{Theorem}\label{Theorem1} \cite[Theorem 5.6]{C3} The positive integer $u=[S/m_S:R/m_R]$ is such that 
$$
\lim_{n\rightarrow \infty}\frac{\ell_R(R/I_n)}{n^d}=u({\rm vol}(\Delta(\hat \Gamma))-{\rm vol}(\Delta(\Gamma))
$$
where 
$$
\Gamma=\{(n_1,\ldots,n_d,i)\in \NN^{d+1}\mid \dim_k \frac{I_i\cap K_{n_1\lambda_1+\cdots+n_d\lambda_d}}{
I_i\cap K^+_{n_1\lambda_1+\cdots+n_d\lambda_d}} >0\mbox{ and }n_1+\cdots+n_d\le\beta i\}
$$
and
$$
\hat\Gamma=\{(n_1,\ldots,n_d,i)\in \NN^{d+1}\mid \dim_k\frac{R\cap K_{n_1\lambda_1+\cdots+n_d\lambda_d}}{
R\cap K^+_{n_1\lambda_1+\cdots+n_d\lambda_d}}>0\mbox{ and }n_1+\cdots+n_d\le\beta i\}.
$$
\end{Theorem}

The sets $\Delta(\Gamma)$ and $\Delta(\hat\Gamma)$ are the closed convex bodies (the Newton-Okounkov bodies) associated  to the semigroups $\Gamma$ and $\hat \Gamma$ as explained in \cite{C1}, \cite{C2} and \cite{C3}.
That is, $\Delta(\Gamma)$ is the intersection of the closed cone in $\RR^{d+1}$ generated by the semigroup $\Gamma$ with $\RR^d\times\{1\}$ and $\Delta(\hat \Gamma)$ is the intersection of the closed cone in $\RR^{d+1}$ generated by the semigroup $\hat \Gamma$ with $\RR^d\times\{1\}$.

By the natural identification of $\RR^d\times\{1\}$ with $\RR^d$,  we will regard $\Delta(\Gamma)$ and $\Delta(\hat\Gamma)$ as convex bodies in $\RR^d$.

\begin{Proposition}\label{Prop1} Suppose that $0\in\Delta(\Gamma)$. Then $\Delta(\Gamma)=\Delta(\hat \Gamma)$.
\end{Proposition}

\begin{proof} We have that $\Delta(\hat\Gamma)$ is the  closure of the set
$$
\left\{\left(\frac{m_1}{i},\ldots,\frac{m_d}{i}\right)\mid (m_1,\ldots,m_d,i)\in \hat\Gamma\right\}
$$
and
$\Delta(\Gamma)$ is the  closure of the set
$$
\left\{\left(\frac{m_1}{i},\ldots,\frac{m_d}{i}\right)\mid (m_1,\ldots,m_d,i)\in \Gamma\right\}.
$$
Since $\Delta(\Gamma)\subset\Delta(\hat\Gamma)$ we must show that if $(m_1,\ldots,m_d,m)\in\hat\Gamma$ and $\epsilon>0$, then there exists $(n_1,\ldots,n_d,n)\in\Gamma$ such that 
$$
\left|\left|\left(\frac{n_1}{n},\ldots,\frac{n_d}{n}\right)-\left(\frac{m_1}{m},\ldots,\frac{m_d}{m}\right)\right|\right|<\epsilon.
$$
Given $(m_1,\ldots,m_d,m)\in \hat\Gamma$, there exists $f\in R$ such that $\nu(f)=m_1\lambda_1+\cdots+m_d\lambda_d$ and $m_1+\cdots+m_d\le\beta  m$.

First suppose that $m_1+\cdots+m_d=\beta m$. Then 
$$
m_1\lambda_1+\cdots+m_d\lambda_d \ge m_1+\cdots +m_d= \beta m
$$
implies $f\in I_m$ by (\ref{eq1}) so $(m_1,\ldots,m_d,m)\in \Gamma$.

Now suppose that $m_1\lambda_1+\cdots+m_d\lambda_d<\beta m$. Since by assumption, $0\in \Delta(\Gamma)$, given $\epsilon>0$, there exists $n>0$ and $g\in I_n$ such that $\nu(g)=n_1\lambda_1+\cdots+n_d\lambda_d$ with $n_1+\cdots+n_d\le \beta n$ and 
$$
\left|\left|\left(\frac{n_1}{n},\ldots,\frac{n_d}{n}\right)\right|\right|<\epsilon.
$$
We can assume that $\epsilon$ is sufficiently small so  that 
\begin{equation}\label{eq2}
\frac{n_1+\cdots+n_d}{n}<\beta-\frac{m_1+\cdots+m_d}{m}.
\end{equation}
We have that $f^ng^m\in I_{mn}$ with
$$
\nu(f^ng^m)=n\nu(f)+m\nu(g)=(nm_1+mn_1)\lambda_1+\cdots+(nm_d+mn_d)\lambda_d.
$$
By (\ref{eq2}), we have that 
$$
\frac{m_1+\cdots+m_d}{m}+\frac{n_1+\cdots+n_d}{n}<\beta
$$
which implies
$$
(nm_1+mn_1)+\cdots+(nm_d+mn_d)=n(m_1+\cdots+m_d)+m(n_1+\cdots+n_d)<mn\beta.
$$
Thus 
$$
\left(\frac{nm_1+mn_1}{mn},\ldots,\frac{nm_d+mn_d}{mn}\right)\in \Delta(\Gamma)
$$
and
$$
\left|\left|\left(\frac{nm_1+mn_1}{mn},\ldots,\frac{nm_d+mn_d}{mn}\right)-\left(\frac{m_1}{m},\ldots,\frac{m_d}{m}\right)\right|\right|
=\left|\left|\left(\frac{n_1}{n},\ldots,\frac{n_d}{n}\right)\right|\right|<\epsilon
$$
so $\left(\frac{m_1}{m},\ldots,\frac{m_d}{m}\right)$ is in the closure of  $\Delta(\Gamma)$ and thus is in $\Delta(\Gamma)$.

\end{proof}

\begin{Lemma}\label{Lemma1} Suppose that 
$$
\lim_{n\rightarrow\infty}\frac{\ell_R(R/I_n)}{n^d}\ne 0.
$$
Then there exists $b\in \ZZ_+$ and $\beta$ as in the equation (\ref{eq1}) such that
\begin{equation}\label{eq3}
I_{ib\beta}\subset m_R^i
\end{equation}
for all $i\in \ZZ_+$.
\end{Lemma}

\begin{proof} By Theorem \ref{Theorem1}, ${\rm vol}(\Delta(\Gamma))<{\rm vol}(\Delta(\hat\Gamma))$ which implies $0\not\in \Delta(\Gamma)$ by Proposition \ref{Prop1}.

Since  $\Delta(\Gamma)$ is closed, there exists $\epsilon>0$ such that the open ball $B_{\epsilon}(0)$ of radius $\epsilon$ centered at 0 in $\RR^d$ is disjoint from $\Delta(\Gamma)$. 

For $c\in \QQ_+$, let $T_c$ be the simplex 
$$
T_c=\{(a_1,\ldots,a_d)\in \RR^d\mid a_i\ge 0\mbox{ for }1\le i\le d\mbox{ and }a_1+\cdots+a_d\le c\}.
$$
Since $B_{\epsilon}(0)\cap \Delta(\Gamma)=\emptyset$, there exists $c\in \QQ_+$ such that $T_c\cap \Delta(\Gamma)=\emptyset$. Thus $n_1+\cdots+n_d>ci$ for all $(n_1,\ldots,n_d,i)\in \Gamma$. We can choose $c$ sufficiently small so that $c<\beta$.

Suppose $f\in I_i$ and $\nu(f)=n_1\lambda_1+\cdots+n_d\lambda_d$. If $n_1+\cdots+n_d<\beta i$, then $(n_1,\ldots,n_d,i)\in \Gamma$ which implies
$$
ic<n_1+\cdots+n_d\le \lambda_1n_1+\cdots+\lambda_dn_d=\nu(f)
$$
so that $f\in K_{ic}\cap R$. If $n_1+\cdots+n_d\ge\beta i$, then $f\in K_{ic}\cap R$ since $\beta\ge c$. Thus
\begin{equation}\label{eq4}
I_i\subset K_{ic}\cap R.
\end{equation}
 Write $c=\frac{a}{b}$ with $a,b\in \ZZ_+$. Then by (\ref{eq1}),
\begin{equation}\label{eq5}
I_{ib\beta}\subset K_{ib\beta c}\cap R=K_{ia\beta}\cap R\subset m_R^{ia\overline c}\subset m_R^i
\end{equation}
for all $i$.
\end{proof}

\section{Positivity of mixed mutiplicities}\label{Sec3}
In this section, we will prove proposition \ref{Prop0}, theorem \ref{Theorem2} and its corollary \ref{Cor3}.  The proof of the general criterion \ref{Theorem3} will be proved in section 5.  We also give an example in this section  to show that the mixed multiplicities of filtrations can be zero even if all the ordinary multiplicities involved are positive in an analytically reducible local ring.

\vskip .2truein
\subsection{Proof of Proposition \ref{Prop0}} Let $\mathcal I_a(j)$ be the $a$-th truncated filtration of $\mathcal I(j)$.
By \cite[Proposition 6.2]{CSS}, 
$$
\lim_{a\rightarrow \infty} e_R(\mathcal I_a(1)^{[d_1]},\ldots,\mathcal I_a(r)^{[d_r]};M)
=e_R(\mathcal I(1)^{[d_1]},\ldots,\mathcal I(r)^{[d_r]};M)
$$
for all $d_1,\ldots,d_r$ with $d_1+\cdots+d_r=d$. Further,  since the $\mathcal I_a(j)$ are Noetherian filtrations, by \cite[Lemma 3.3]{CSS} each $e_R(\mathcal I_a(1)^{[d_1]},\ldots,\mathcal I_a(r)^{[d_r]};M)$ is a positive constant times a mixed multiplicity $e_R( J_1^{[d_1]},\ldots, J_r^{[d_r]};M)$ of a set of $m_R$-primary ideals $J_1,\ldots,J_r$ (which depend on $a$). This mixed multiplicity is nonnegative by \cite{T1} or  \cite[Corollary 17.4.7]{HS}.

\subsection{Proof of Theorem \ref{Theorem2}}

 Since  
 $$
 \ell_R(R/I(1)_{n_1}\cdots I(1)_{n_r})=\ell_{\hat R}(\hat R/I(1)_{n_1}\cdots I(1)_{n_r}\hat R) 
 $$
 for $n_1,\ldots,n_r\in \NN$, we may assume that $R$ is a complete domain. 
 
  By  \cite[Lemma 3.3]{CSS},     we have equality of  mixed multiplicities 
   \begin{equation}\label{eq6}
 e_R(\mathcal I(1)_t^{[d_1]},\ldots,\mathcal I(r)_t^{[d_r]};R)=\frac{1}{s_t^d}e_R((I(1)_{t,s_t})^{[d_1]},\ldots,(I(r)_{t,s_t})^{[d_r]};R)
 \end{equation}
  where $s_t$ is such that $I(j)_{t,s_ti}=(I(j)_{t,s_t})^i$ for all $i>0$ and $1\le j\le r$. By (\ref{eq3}), there exists $b\in \NN$ such that $I(j)_{ib\beta}\subset m_R^i$ for $1\le j\le r$ and $i\in \ZZ_+$.  Thus
 \begin{equation}
 I(j)_{t,s_ti}\subset I(j)_{s_ti}\subset m_R^{i\frac{s_t}{b\beta}}\mbox{ for all $i$ and $j$}
 \end{equation}
 if $s_t$ is chosen to be a multiple of $b\beta$. Thus 
 $$
 \begin{array}{lll}
 e_R((I(1)_{t,s_t})^{[d_1]},\ldots,(I(r)_{t,s_t})^{[d_r]};R)&\ge& e_R((m_R)^{\frac{s_t}{b\beta}})^{[d_1]},\ldots,(m_R^{\frac{s_t}{b\beta}})^{[d_r]};R)\\
 &=&\frac{(s_t)^d}{(b\beta)^d}e_R(m_R^{[d_1]},\ldots,m_R^{[d_r]};R)
 \end{array}
 $$
for all $t$ and $d_1,\ldots,d_r\in \NN$ with $d_1+\cdots+d_r=d$ by the inequality of mixed multiplicities of $m_R$-primary ideals of  \cite[Lemma 17.5.3]{HS} or \cite[Lemma 14, page 8]{KG}, so 
$$
e_R(\mathcal I(1)_t^{[d_1]},\ldots,\mathcal I(r)_t^{[d_r]};R)\ge \frac{1}{(b\beta)^d}e_R(m_R^{[d_1]},\ldots,m_R^{[d_r]};R)
$$
for all $t$, $d_1,\ldots,d_r\in \NN$ with $d_1+\cdots+d_r=d$ by (\ref{eq6}). Thus
$$
e_R(\mathcal I(1)^{[d_1]},\ldots,\mathcal I(r)^{[d_r]};R)=\lim_{t\rightarrow\infty}e_R(\mathcal I(1)_t^{[d_1]},\ldots,\mathcal I(r)_t^{[d_r]};R)
\ge \frac{1}{(b\beta)^d}e_R(m_R^{[d_1]},\ldots,m_R^{[d_r]};R)
$$
for all $d_1,\ldots,d_r$ with $d_1+\cdots+d_r=d$ by \cite[Proposition 6.2]{CSS}. Finally, we observe that each mixed multiplicity 
$e_R(m_R^{[d_1]},\ldots,m_R^{[d_r]};R)$ is the ordinary multiplicity $e_R(m_R;R)$ of $R$, and hence is positive.

\subsection{Proof of Corollary \ref{Cor3}}

 By \cite[Theorem 6.8]{CSS},
for any $d_1,\ldots,d_r\in \NN$ with $d_1+\cdots+d_r=d$, 
$$
\begin{array}{lll}
e_R(\mathcal I(1)^{[d_1]},\ldots,\mathcal I(r)^{[d_r]};M)&=&
e_{\hat R}((\mathcal I(1)\hat R)^{[d_1]},\ldots,(\mathcal I(r)\hat R)^{[d_r]};\hat M)\\
&=&
\sum\ell_{\hat R_P}(\hat M_P)e_{\hat R/P}((\mathcal I(1)\hat R/P)^{[d_1]},\ldots,(\mathcal I(r)\hat R/P)^{[d_r]};\hat R/P)
\end{array}
$$
where the sum is  over the minimal primes $P$ of $\hat R$ such that $\dim \hat R/P=d$ and $\mathcal I(j)\hat R/P=\{I(j)_i\hat R/P\}$. The corollary now follows from Theorem \ref{Theorem2}.

\subsection{Construction  of an example}
\begin{Example}\label{Example1} There exists a two-dimensional excellent local domain $R$ and filtrations $\mathcal I$ and $\mathcal J$ of $R$ by  $m_R$-primary ideals such that $e_R(\mathcal I^{[2]};R)=e_R(\mathcal I;R)>0$, 
$e_R(\mathcal J^{[2]};R)=e_R(\mathcal J;R)>0$, but the mixed multiplicity $e_R(\mathcal I^{[1]},\mathcal J^{[1]};R)=0$.
\end{Example} 

\begin{proof}
 Let $R=\CC[x,y,z]_{(x,y,z)}/(y^2-x^2(x+1))$, which is a two dimensional excellent domain. The minimal primes of the $m_R$-adic completion of $R$ are $P_1=(y-x\sqrt{x+1})$ and $P_2=(y+x\sqrt{x+1})$. Let $R_1=\hat R/(y-x\sqrt{x+1})\cong \CC[[x,z]]$ and $R_2=\hat R/(y+x\sqrt{x+1})\cong \CC[[x,z]]$. By \cite[Lemma 5.1]{C1}, if $\mathcal I=\{I_n\}$ is a graded family of $m_R$-primary ideals, then 
\begin{equation}\label{eqL}
\lim_{n\rightarrow \infty}\frac{\ell_R(R/I_n)}{n^2}
=\lim_{n\rightarrow \infty}\frac{\ell_{\hat R}(\hat R/I_n\hat R)}{n^2}
=\lim_{n\rightarrow\infty}\frac{\ell_{R_1}(R_1/I_nR_1)}{n^2}+
\lim_{n\rightarrow\infty}\frac{\ell_{R_2}(R_2/I_nR_2)}{n^2}.
\end{equation}
We have the expansion
$$
x\sqrt{x+1}=a_1x+a_2x^2+a_3x^3+\cdots
$$
where 
$$
a_{n+1}=\frac{(-1)^{n-1}(2n-3)!}{2^{2n-2}n!(n-2)!}.
$$
Define filtrations of $m_R$-primary ideals by $\mathcal I=\{I_n\}$ with 
$$
I_n=(y-a_1x-a_2x^2-\cdots-a_{n-1}x^{n-1})+m_R^n
$$
and $\mathcal J=\{J_n\}$ with 
$$
J_n=(y+a_1x+a_2x^2+\cdots+a_{n-1}x^{n-1})+m_R^n.
$$
We have that
$$
I_nR_1=m_{R_1}^n, I_nR_2=(x)+m_{R_2}^n, J_nR_1=(x)+m_{R_1}^n, J_nR_2=m_{R_2}^n
$$
for $n\ge 1$, so that 
$$
I_nJ_nR_1=xm_{R_1}^n+m_{R_1}^{2n}.
$$
The set of all monomials in $x^iz^j$ with $i+j\le n$ and the $n-1$ monomials $z^{n+1},\ldots,z^{2n-1}$ is thus a  $\CC$-basis of  $R_1/I_nJ_nR_1$.  
Further,
$$
I_nJ_nR_2=xm_{R_2}^n+m_{R_2}^{2n},
$$
so the set of all monomials in $x^iz^j$ with $i+j\le n$ and the $n-1$ monomials $z^{n+1},\ldots,z^{2n-1}$ is also a  $\CC$-basis of $R_2/I_nJ_nR_2$.  

Thus
$$
\lim_{n\rightarrow \infty}\frac{\ell_{R_1}(R_1/I_nR_1)}{n^2}=\frac{1}{2},
\lim_{n\rightarrow \infty}\frac{\ell_{R_2}(R_2/I_nR_2)}{n^2}=0,
$$
$$
\lim_{n\rightarrow \infty}\frac{\ell_{R_1}(R_1/J_nR_1)}{n^2}=0,
\lim_{n\rightarrow \infty}\frac{\ell_{R_2}(R_2/J_nR_2)}{n^2}=\frac{1}{2},
$$
$$
\lim_{n\rightarrow \infty}\frac{\ell_{R_1}(R_1/I_nJ_nR_1)}{n^2}=\frac{1}{2},
\lim_{n\rightarrow \infty}\frac{\ell_{R_2}(R_2/I_nJ_nR_2)}{n^2}=\frac{1}{2}.
$$
Thus by (\ref{eqL}), $e_R(\mathcal I^{[2]};R)=e_R(\mathcal I;R)=1$ and 
$e_R(\mathcal J^{[2]};R)=e_R(\mathcal J;R)=1$.

Further, we have by (\ref{eqL}) that $\lim_{n\rightarrow\infty}\frac{\ell_R(R/I_nJ_n)}{n^2}=1$.
Now, from \cite[Theorem 6.6]{CSS}, we calculate
$$
\begin{array}{lll}
1=\lim_{n\rightarrow\infty}\frac{\ell_R(R/I_nJ_n)}{n^2}&=&\frac{e_R(\mathcal I^{[2]};R)}{2}+e_R(\mathcal I^{[1]},\mathcal J^{[1]};R)+\frac{e_R(\mathcal J^{[2]};R)}{2}\\
&=& 1+e_R(\mathcal I^{[1]},\mathcal J^{[1]};R)
\end{array}
$$
and conclude that $e_R(\mathcal I^{[1]},\mathcal J^{[1]};R)=0$.
\end{proof}

\section{Minkowski sums of Okounkov bodies} We continue in this section with the notation of Section \ref{Sec2}. In particular, we assume that $R$ is a complete Noetherian local domain. Let $\mathcal I(1)=\{I(1)_n\},\ldots,\mathcal I(r)=\{I(r)_n\}$ be filtrations of $R$ by $m_R$-primary ideals. 
For all $(\sigma_1,\ldots,\sigma_r)\in \NN^r$, define  semigroups
\begin{equation}\label{eqX8}
\begin{array}{lll}
\Gamma_{(\sigma_1,\ldots,\sigma_r)}
&=&\{(n_1,\ldots,n_d,i)\in \NN^{d+1}\mid \\
&&\dim_kI(1)_{i\sigma_1}\cdots I(r)_{i\sigma_r}\cap K_{n_1\lambda_1+\cdots+n_d\lambda_d}
/I(1)_{i\sigma_1}\cdots I(r)_{i\sigma_r}\cap K^+_{n_1\lambda_1+\cdots+n_d\lambda_d}>0\\
&& \mbox{ and }n_1+\cdots +n_d\le\beta i\},
\end{array}
\end{equation}
where $\beta$ is chosen so that (\ref{eq1}) holds for $I_n=I(a)_{n\sigma_1}\cdots I(r)_{n\sigma_r}$.
With the notation of Section \ref{Sec2}, we have that $\hat \Gamma=\Gamma_{(0,\ldots,0)}$.

\begin{Lemma}\label{LemmaX2} 
Suppose $(\sigma_1,\ldots,\sigma_r),(\tau_1,\ldots,\tau_r)\in\NN^r$ are such that
$$
\Delta(\Gamma_{(\sigma_1,\ldots,\sigma_r)})\subset \Delta(\Gamma_{(\tau_1,\ldots,\tau_r)})
$$
(with $\beta \gg 0$ in (\ref{eqX8})) and
$$
{\rm vol}(\Delta(\Gamma_{(\sigma_1,\ldots,\sigma_r)}))={\rm vol}(\Delta(\Gamma_{(\tau_1,\ldots,\tau_r)})).
$$
Then 
$$
\Delta(\Gamma_{(\sigma_1,\ldots,\sigma_r)})= \Delta(\Gamma_{(\tau_1,\ldots,\tau_r)}).
$$
\end{Lemma}

\begin{proof}
Suppose $\Delta(\Gamma_{(\tau_1,\ldots,\tau_r)})\ne \Delta(\Gamma_{(\sigma_1,\ldots,\sigma_r)})$.  Then there exists 
  $$
  p\in \Delta(\Gamma_{(\tau_1,\ldots,\tau_r)}) \setminus
  \Delta(\Gamma_{(\sigma_1,\ldots,\sigma_r)}).
  $$
  Since $\Delta(\Gamma_{(\sigma_1,\ldots,\sigma_r)})$ is closed in $\RR^d$, there exists an epsilon ball $B_{\epsilon}(p)$ centered at $p$ in $\RR^d$ such that $B_{\epsilon}(p)\cap \Delta(\Gamma_{(\sigma_1,\ldots,\sigma_r)})=\emptyset$. Now $\Delta(\Gamma_{(\tau_1,\ldots,\tau_r)})$ has positive volume (since $\beta\gg0$) so there exist $w_1,\ldots,w_d\in \Delta(\Gamma_{(\tau_1,\ldots,\tau_r)})$ such that $v_1=w_1-p,\ldots,v_d=w_d-p$ is a real basis of $\RR^d$. Since $\Delta(\Gamma_{(\tau_1,\ldots,\tau_r)})$ is convex, there exists $\delta>0$ such that letting $W$ be the hypercube 
  $$
  W=\{p+\alpha_1 v_1+\cdots +\alpha_dv_d\mid 0\le \alpha_i\le \delta\mbox{ for }1\le i\le d\},
  $$
  we have that 
  $$
  W\subset \Delta(\Gamma_{(\tau_1,\ldots,\tau_r)})\cap B_{\epsilon}(p).
  $$
  But then 
  $$
 {\rm vol}( \Delta(\Gamma_{(\tau_1,\ldots,\tau_r)}) )
 -{\rm vol}( \Delta(\Gamma_{(\sigma_1,\ldots,\sigma_r)}))\ge {\rm vol}(W)>0,
 $$
 a contradiction. Thus
 $$
  \Delta(\Gamma_{(\sigma_1,\ldots,\sigma_r)}) 
 =\Delta(\Gamma_{(\tau_1,\ldots,\tau_r)}).
 $$
 \end{proof}

Let ${\rm HF}$ be the half space
$$
{\rm HF}=\{(a_1,\ldots,a_d)\in \RR^d\mid a_1+\cdots+a_d\le\beta\}.
$$
\begin{Lemma}\label{LemmaX1} 
For $(\sigma_1,\ldots,\sigma_r), (\tau_1,\ldots,\tau_r)\in \NN^r$ (with $\beta\gg 0$ in (\ref{eqX8})) we have that
$$
\left[\Delta(\Gamma_{(\sigma_1,\ldots,\sigma_r)})+\Delta(\Gamma_{(\tau_1,\ldots,\tau_r)})\right]\cap {\rm HF}
\subset \Delta(\Gamma_{(\sigma_1+\tau_1,\ldots,\sigma_r+\tau_r)})
$$
where $\Delta(\Gamma_{(\sigma_1,\ldots,\sigma_r)})+\Delta(\Gamma_{(\tau_1,\ldots,\tau_r)})$ is the Minkowski sum of $\Delta(\Gamma_{(\sigma_1,\ldots,\sigma_r)})$ and $\Delta(\Gamma_{(\tau_1,\ldots,\tau_r)})$.
\end{Lemma}

\begin{proof} The set $\left[\Delta(\Gamma_{(\sigma_1,\ldots,\sigma_r)})+\Delta(\Gamma_{(\tau_1,\ldots,\tau_r)})\right]\cap {\rm HF}$ is the closure of the set of points 
$$
\left(\frac{m_1}{i},\ldots,\frac{m_d}{i}\right)+\left(\frac{n_1}{j},\ldots,\frac{n_d}{j}\right)
$$
such that $(m_1,\ldots,m_d,i)\in\Gamma_{(\sigma_1,\ldots,\sigma_r)}$, 
$(n_1,\ldots,n_d,j)\in\Gamma_{(\tau_1,\ldots,\tau_r)}$ and
\begin{equation}\label{eqX7}
\frac{jm_1+in_1}{ij}+\cdots+\frac{jm_d+in_d}{ij}\le\beta.
\end{equation}
It thus  suffices to show that if $(m_1,\ldots,m_d,i)\in \Gamma_{(\sigma_1,\ldots,\sigma_r)}$ and $(n_1,\ldots,n_d,j)\in \Gamma_{(\tau_1,\ldots,\tau_r)}$ satisfy (\ref{eqX7}), then
$$
\left(\frac{m_1}{i},\ldots,\frac{m_d}{i}\right)+\left(\frac{n_1}{j},\ldots,\frac{n_d}{j}\right)\in \Delta(\Gamma_{(\sigma_1+\tau_1,\ldots,\sigma_r+\tau_r)}).
$$
Assume $(m_1,\ldots,m_d,i)$ and $(n_1,\ldots,n_d,j)$ satisfy these conditions. Then there exists $f\in I(1)_{i\sigma_1}\cdots I(r)_{i\sigma_r}$ such that $\nu(f)=m_1\lambda_1+\cdots+m_d\lambda_d$ with $m_1+\cdots+m_d\le i\beta$ and there exists $g\in I(1)_{j\tau_1}\cdots I(r)_{j\tau_r}$ such that 
$\nu(g)=n_1\lambda_1+\cdots+n_d\lambda_d$ with $n_1+\cdots+n_d\le j\beta$. Then
$$
f^jg^i\in I(1)_{ji(\sigma_1+\tau_1)}\cdots I(r)_{ji(\sigma_r+\tau_r)}
$$
with $\nu(f^jg^i)=\lambda_1(jm_1+in_1)+\cdots+\lambda_d(jm_d+i n_d)$ and 
$(jm_1+in_1)+\cdots+(jm_d+in_d)\le ij\beta$ by (\ref{eqX7}). Thus
$$
\left(\frac{m_1}{i},\ldots,\frac{m_d}{i}\right)+\left(\frac{n_1}{j},\ldots,\frac{n_d}{j}\right)\in \Delta(\Gamma_{(\sigma_1+\tau_1,\ldots,\sigma_r+\tau_r)}).
$$
\end{proof}

\begin{Proposition}\label{PropX2} Suppose that $(\sigma_1,\ldots,\sigma_r),(\tau_1,\ldots,\tau_r)\in \NN^r$ (with $\beta\gg 0$ in (\ref{eqX8}))
and $\Delta(\Gamma_{(\sigma_1,\ldots,\sigma_r)})=\Delta(\Gamma_{(0,\ldots,0)})$. Then
$\Delta(\Gamma_{(\sigma_1+\tau_1,\ldots,\sigma_r+\tau_r)})=\Delta(\Gamma_{(\tau_1,\ldots,\tau_r)})$.
\end{Proposition}

\begin{proof} By Lemma \ref{LemmaX1},
$$
\left[\Delta(\Gamma_{(0,\ldots,0)})+\Delta(\Gamma_{(\tau_1,\ldots,\tau_r)})\right]\cap {\rm HF}
\subset \Delta(\Gamma_{(\tau_1,\ldots,\tau_r)}).
$$
Now since $0\in\Delta(\Gamma_{(0,\ldots,0)})$, we have
$$
\Delta(\Gamma_{(\tau_1,\ldots,\tau_r)})\subset \left[\Delta(\Gamma_{(0,\ldots,0)})+\Delta(\Gamma_{(\tau_1,\ldots,\tau_r)})\right]\cap {\rm HF}.
$$
Thus
$$
\Delta(\Gamma_{(\tau_1,\ldots,\tau_r)})\subset \left[\Delta(\Gamma_{(\sigma_1,\ldots,\sigma_r)})+\Delta(\Gamma_{(\tau_1,\ldots,\tau_r)})\right]\cap {\rm HF}
$$
and so
\begin{equation}\label{eqX1}
\Delta(\Gamma_{(\tau_1,\ldots,\tau_r)})\subset\Delta(\Gamma_{(\sigma_1+\tau_1,\ldots,\sigma_r+\tau_r)})
\end{equation}
 by Lemma \ref{LemmaX1}. Thus
 $$
 {\rm vol}(\Delta(\Gamma_{(\tau_1,\ldots,\tau_r)}))\le {\rm vol}(\Delta(\Gamma_{(\sigma_1+\tau_1,\ldots,\sigma_r+\tau_r)}))
 $$
and so
$$
{\rm vol}(\Delta(\Gamma_{(0,\ldots,0)}))-{\rm vol}(\Delta(\Gamma_{(\sigma_1+\tau_1,\ldots,\sigma_r+\tau_r)}))
\le 
{\rm vol}(\Delta(\Gamma_{(0,\ldots,0)}))-{\rm vol}(\Delta(\Gamma_{(\tau_1,\ldots,\tau_r)})).
$$
Thus by Theorem \ref{Theorem1},
\begin{equation}\label{eqX4}
\lim_{t\rightarrow\infty}\frac{\ell_R(R/ I(1)_{t(\sigma_1+\tau_1)}\cdots I(r)_{t(\sigma_r+\tau_r)})}{t^d}
\le 
\lim_{t\rightarrow\infty}\frac{\ell_R(R/ I(1)_{t\tau_1}\cdots I(r)_{t\tau_r})}{t^d}.
 \end{equation}
 Now for all $t\in \ZZ_+$, there are natural surjections
 $$
 R/I(1)_{t(\sigma_1+\tau_1)}\cdots I(r)_{t(\sigma_r+\tau_r)}\rightarrow R/I(1)_{t\tau_1}\cdots I(r)_{t\tau_r}
 $$
 which implies 
 \begin{equation}\label{eqX3}
 \lim_{t\rightarrow\infty}\frac{\ell_R(R/ I(1)_{t\tau_1}\cdots I(r)_{t\tau_r})}{t^d}
 \le
 \lim_{t\rightarrow\infty}\frac{\ell_R(R/ I(1)_{t(\sigma_1+\tau_1)}\cdots I(r)_{t(\sigma_r+\tau_r)})}{t^d}.
  \end{equation}
  Thus
  $$
  \lim_{t\rightarrow\infty}\frac{\ell_R(R/ I(1)_{t(\sigma_1+\tau_1)}\cdots I(r)_{t(\sigma_r+\tau_r)})}{t^d} 
 =\lim_{t\rightarrow\infty}\frac{\ell_R(R/ I(1)_{t\tau_1}\cdots I(r)_{t\tau_r})}{t^d}
 $$
  by (\ref{eqX4}) and  (\ref{eqX3}). By Theorem \ref{Theorem1}, we have that 
  $$
  {\rm vol}(\Delta(\Gamma_{(\sigma_1+\tau_1,\ldots,\sigma_r+\tau_r)}))=
  {\rm vol}(\Delta(\Gamma_{(\tau_1,\ldots,\tau_r)})),
  $$
  and so
   $$
  \Delta(\Gamma_{(\sigma_1+\tau_1,\ldots,\sigma_r+\tau_r)}) 
 =\Delta(\Gamma_{(\tau_1,\ldots,\tau_r)})
 $$
 by (\ref{eqX1}) and Lemma \ref{LemmaX2}.
 \end{proof}
 
 \section{Proof of Theorem \ref{Theorem3}}\label{Sec5}
 
 Since $\ell_R(M/I(1)_{n_1}\cdots I(r)_{n_r}M)=\ell_{\hat R}(\hat M/I(1)_{n_1}\cdots I(r)_{n_r}\hat M) $ for all $n_1,\ldots,n_r\in \NN^r$, we may assume that $R$ is a complete domain. By \cite[Theorem 6.8]{CSS} we may assume that $M=R$.
 
 The assumption $e_R(\mathcal I(j);R)=0$ for $j>s$ implies (taking $\beta\gg0$ in (\ref{eqX8})) by Theorem \ref{Theorem1} and Lemma \ref{LemmaX2} that 
 $$
 \Delta(\Gamma_{(0,\ldots,0,1,0,\ldots,0)})=\Delta(\Gamma_{(0,\ldots,0)})
 $$
 whenever the 1 is in a position greater than $s$.
 
 By Proposition \ref{PropX2} and Theorem \ref{Theorem1}, we have that
 \begin{equation}\label{eqX6}
 \lim_{t\rightarrow\infty}\frac{\ell_R(R/I(1)_{tn_1}\cdots I(r)_{tn_r})}{t^d}
 = \lim_{t\rightarrow\infty}\frac{\ell_R(R/I(1)_{tn_1}\cdots I(s)_{tn_s})}{t^d}
 \end{equation}
 for all $(n_1,\ldots,n_r)\in \NN^r$. The function 
 $$
 G(n_1,\ldots,n_r)=\lim_{t\rightarrow \infty}\frac{\ell_R(R/I(1)_{tn_1}\cdots I(r)_{tn_r})}{t^d}
 $$
 is a homogeneous polynomial in $n_1,\ldots,n_r$ of total degree $d$ by \cite[Theorem 6.6]{CSS} (recalled in (\ref{M2}) of this paper). The mixed multiplicities are defined from this polynomial by the writing
 \begin{equation}\label{eqV6}
 G(n_1,\ldots,n_r)=\sum_{d_1+\cdots+d_r=d}\frac{1}{d_1!\cdots d_r!}
 e_R(\mathcal I(1)^{[d_1]},\ldots,\mathcal I(r)^{[d_r]};R)n_1^{d_1}\cdots n_r^{d_r}.
 \end{equation}
 By (\ref{eqX6}), we have that $G$ does not depend on $n_{s+1},\ldots,n_r$ so that (\ref{eqX41}) holds.

Equation (\ref{eqX5}) follows from \cite[Proposition 6.5]{CSS} and Theorem \ref{Theorem2}.


\end{document}